\documentclass{article}

\usepackage{amsmath,amssymb}
\usepackage{enumerate}
\usepackage{url}
\usepackage{comment}
\usepackage{theorem} 
\theoremstyle{break} \newtheorem{Theorem}{Theorem} 
\theoremstyle{break} \newtheorem{Definition}[Theorem]{Definition} 
\theoremstyle{break} \newtheorem{Lemma}[Theorem]{Lemma} 
\theoremstyle{break} \newtheorem{Corollary}[Theorem]{Corollary} 
\theoremstyle{break} \newtheorem{Result}[Theorem]{Result} 
\theoremstyle{break} \newtheorem{Example}{Example} 
\theoremstyle{break} \newtheorem{Construction}{Construction}

\newenvironment{Proof}{{\noindent \bf Proof.}}{\par\smallbreak}
\newcommand{\qed}{\hfill\ensuremath{\Box}} 
 \title{Generalised Veroneseans}

\author{A. Klein, J. Schillewaert and L. Storme }
\begin{document}
\maketitle
\begin{abstract}
In \cite{ThasHVM}, a characterization of the finite quadric Veronesean
$\mathcal{V}_{n}^{2^{n}}$ 
by means of properties of the set of its tangent spaces is
proved. These tangent spaces form a {\em regular generalised dual
  arc}. 
We prove an extension result for regular generalised dual arcs.
To motivate our research, we show how they  are used to construct a large class of secret sharing schemes.  
\end{abstract}  

A typical problem in (finite) geometry is the study of highly symmetrical substructures. For example, arcs are configurations of points in PG$(n,q)$ such that each $n+1$ of them are in general position, while $n_{1}$-dimensional dual arcs are sets of $n_{1}$-spaces such that each two intersect in a point and any three of them are skew.  
These two structures appear naturally in cryptographical applications.

In this article, we define objects, called (generalised) dual arcs; a class of objects that contain classical arcs and $n_{1}$-dimensional dual arcs as special cases. These (generalised) dual arcs have applications in cryptography as well. 

We give a construction method for a wide class of parameters and prove an extension result for regular generalised dual arcs of order $d=1$.

In Sections \ref{DC} and \ref{ss}, we give the necessary definitions, constructions, and examples of applications in cryptography. Section \ref{kr} refers to known classification results, and Section \ref{sec:d1} states our main characterization theorem (Theorem \ref{T:d1}). 
We now start with the required definitions to make this article self-contained.

\section{Definitions and  constructions} \label{DC}

\begin{Definition}
  A {\em generalised dual arc $\mathcal{F}$} of order $d$ with dimensions
  $n=n_{0}> n_1 > n_2 > \cdots > n_{d+1}>-1$ of $PG(n,q)$ is a set of $n_{1}$-dimensional subspaces of
  $PG(n,q)$ such that:
  \begin{enumerate}
  \item each $j$ of these subspaces intersect in a subspace of
    dimension $n_j$, $1 \leq j \leq d+1$,
  \item each $d+2$ of these subspaces have no common intersection.

  \end{enumerate}

  We call $(n=n_{0},n_{1},\ldots,n_{d+1})$ the {\em parameters} of the
  generalised dual arc.
\end{Definition}

\begin{Definition}
  A generalised dual arc of order $d$ with parameters
  $(n=n_{0},\ldots,n_{d+1})$ is \emph{regular} if, in addition, the
  $n_1$-dimensional spaces span $PG(n,q)$ and if it satisfies 
  the property that if $\pi$ is the intersection of $j$ elements of
  $\mathcal{F}$, $j \leq d$, then $\pi$ is spanned by the subspaces of
  dimension $n_{j+1}$ which are the intersections of $\pi$ with the
  remaining elements of $\mathcal{F}$.
\end{Definition}

\begin{Construction}\label{cons:1}
  Let $PG(V)$ be an $n$-dimensional space with basis $e_i$ ($0 \leq i \leq
  n$).

  Let $PG(W)$ be an $\left(\binom{n+d+1}{d+1}-1\right)$-dimensional space
  with basis $e_{i_0,\ldots,i_d}$ ($0 \leq i_0 \leq i_1 \leq \cdots \leq i_d \leq n$).
  
  We now define a multilinear mapping from $PG(V)$ to $PG(W)$. In the description of this multilinear mapping and in the remainder of this article, the vector 
   $e_{i_0,\ldots,i_d}$, for $0 \leq i_0 , i_1, \ldots , i_d \leq n$, is identical to the vector  $e_{i_{\sigma(0)},  \ldots, i_{\sigma(d)}}$, where $\sigma$ is a permutation of $\{0,\ldots,d\}$ with $0\leq i_{\sigma(0)} \leq  \cdots \leq i_{\sigma(d)}\leq n$. For example, $e_{001},e_{010},e_{001}$ all denote the same vector  $e_{001}$.

  Let $\theta: V^{d+1} \to W$ be the multilinear mapping
  \begin{equation}
    \label{eq:theta}
    \theta: (\sum_{i_{0}=0}^{n}x_{i_{0}}^{(0)}e_{i_{0}},\ldots,\sum_{i_{d}=0}^{n}x_{i_{d}}^{(d)}e_{i_{d}}) 
    \mapsto \sum_{0\leq i_{0},\ldots ,i_{d} \leq n}
    x_{i_0}^{(0)}\cdot\ldots\cdot x_{i_d}^{(d)} e_{i_0,\ldots,i_d} \ .
  \end{equation}

For example, $\theta(x_0^{(0)}e_0+x_1^{(0)}e_1,x_0^{(1)}e_0+x_1^{(1)}e_1)=x_0^{(0)}x_0^{(1)}e_{0,0}+(x_0^{(0)}x_1^{(1)}+x_1^{(0)}x_0^{(1)})e_{0,1}+x_1^{(0)}x_1^{(1)}e_{1,1}$.

  For each point $P=[x]$ of $PG(V)$, we define a subspace $D(P)$ of
  $PG(W)$ by
  \begin{equation}
    \label{eq:Dtheta}
    D(P) = \left< \theta(x,v_1,\ldots,v_d)\mid v_1,\ldots,v_d \in V\right> \ .
  \end{equation}.
\end{Construction}
  
\begin{Theorem}[\cite{KLJ}] \label{T:dual-arc} The set $\mathcal{D} = \{D(P) \mid P \in
  PG(V)\}$ is a generalised dual arc with dimensions $d_i =
  \binom{n+d+1-i}{d+1-i}-1, i=0,\ldots ,d+1$.
\end{Theorem}


For $q$ odd and $\frac{q^{n}-1}{q-1} \geq \binom{n+d}{d+1}$, there is an alternative
construction.

\begin{Construction} \label{C:2}
  We define $\zeta:PG(V) \to PG(W)$ by
  $$ \zeta: [\sum_{i=0}^{n} x_ie_i] \mapsto 
  [\sum_{0 \leq i_0 \leq \cdots \leq i_d \leq n} x_{i_0}\cdots x_{i_d} e_{i_0,\ldots,i_d}].$$

  This mapping $\zeta$ is a generalisation of the well-known
  quadratic Veronesean map (see \cite{MR1363259}). We call it 
  the \emph{generalised Veronesean}.
  
  With $b$ and $B$ respectively, we denote the standard scalar product
  of $V$ and $W$, i.e.,
  $$ b(\sum_{i=0}^{n} x_i e_i,\sum_{i=0}^{n} y_i e_i)=\sum_{i=0}^{n} x_i y_i,$$  
  and
  $$ B(\sum_{0 \leq i_0 \leq \cdots \leq i_d \leq n} x_{i_0,\ldots,i_d}
  e_{i_0,\ldots,i_d},\sum_{0 \leq i_0 \leq \cdots \leq i_d \leq n} y_{i_0,\ldots,i_d} e_{i_0,\ldots,i_d}) =$$
  $$\sum_{0 \leq i_0 \leq \cdots \leq i_d \leq n} x_{i_0,\ldots,i_d}y_{i_0,\ldots,i_d}$$
  with summations over $0\leq i_{0}\leq\cdots\leq i_{d}\leq n$.

  For each $x \in V$, we denote by $x^{\perp}$ the subspace of $V$
  perpendicular to $x$ with respect to $b$. So
  $$ x^{\perp} = \{ y\in V \mid b(x,y)=0\}.$$

  Then 
  \begin{equation}
    \label{eq:D(P)}
    D(P) = \{[z] \in PG(W) \mid B(z,\zeta(y))=0 \text{ for all } y \in x^{\perp}\}
  \end{equation}
  is a generalised dual arc.
\end{Construction}

We will not use Construction~\ref{C:2} in this article. Hence, we
refer to \cite{KLJ} for a proof that Construction~\ref{C:2} gives 
generalised dual arcs isomorphic to the ones described by
Construction~\ref{cons:1}.  

For the second construction, we call the arcs of form $\mathcal{D} =
\{D(P)|P \in PG(V)\}$ \emph{Veronesean dual arcs}.

Below, we give two examples of our general construction.

\begin{Example} \label{E:1} Starting with $PG(2,q)$, the mapping
  $\zeta:PG(2,q) \to PG(5,q)$ with

$$\zeta([x_{0},x_{1},x_{2}])=[x_{0}^{2},x_{1}^{2},x_{2}^{2},x_{0}x_{1},x_{0}x_{2},x_{1}x_{2}]$$
defines the quadratic Veronesean $\mathcal{V}_2^4$.

If $P=[a,b,c]$, the planes $D(P)$ defined above have the equation
   $$ D(P) =
   \{[ax_{0},bx_{1},cx_{2},ax_{1}+bx_{0},ax_{2}+cx_{0},bx_{2}+cx_{1}] \mid
   x_{0},x_{1},x_{2} \in \mathbb{F}_{q}\} \ .$$

   These planes form a regular generalised dual arc of
   $q^2+q+1$ planes with parameters $(5,2,0)$. 
 \end{Example}

 \begin{Example} \label{E:2} The map $\zeta:PG(2,q) \to PG(9,q)$ with

$$\zeta([x_{0},x_{1},x_{2}])=[x_{0}^{3},x_{1}^{3},x_{2}^{3},x_{0}^{2}x_{1},x_{0}^{2}x_{2},x_{1}^{2}x_{0},x_{1}^{2}x_{2},
x_{2}^{2}x_{0},x_{2}^{2}x_{1},x_{0}x_{1}x_{2}]$$ defines a cubic
Veronesean. Construction~\ref{cons:1} associates to each of the
$q^{2}+q+1$ points a $5$-dimensional space in $PG(9,q)$. Each
two of these $5$-spaces intersect in a plane. Each three $5$-spaces
share a common point and each four $5$-spaces are skew.

\end{Example}

Three of the $q^{2}+q+1$ 5-spaces are:
\begin{align*}
  \pi_{0}:=D([1,0,0]) &= \{[e_{0},0,0,e_{1},e_{2},e_{3},0,e_{4},0,e_{5}]
  \mid
  e_{i} \in \mathbb{F}_{q}\}, \\
  \pi_{1}:=D([0,1,0]) &= \{[0,e_{0},0,e_{1},0,e_{2},e_{3},0,e_{4},e_{5}]
  \mid
  e_{i} \in \mathbb{F}_{q}\}, \\
  \pi_{2}:=D([0,0,1]) &= \{[0,0,e_{0},0,e_{1},0,e_{2},e_{3},e_{4},e_{5}]
  \mid
  e_{i} \in \mathbb{F}_{q}\} .\\
\end{align*}

In each $5$-space, the other $q^{2}+q$ $5$-spaces intersect in a
configuration of $q^{2}+q$ planes. These planes are a part of the
Veronesean described in Example \ref{E:1}.

For $\pi_{0}$, the corresponding Veronesean has the form
  $$\mathcal{V}_{0} := [x_{0}^{2},0,0,x_{0}x_{1},x_{0}x_{2},x_{1}^{2},0,x_{2}^{2},0,x_{1}x_{2}].$$

  This Veronesean $\mathcal{V}_0$ has $q^{2}+q+1$ tangent planes; where $q^{2}+q$ of the
  tangent planes are intersections of $\pi_{0}$ with the other $5$-spaces.
  The extra plane has the form
  $$E_{0}:=\{[e_{0},0,0,e_{1},e_{2},0,0,0,0,0] \mid e_{0},e_{1},e_{2} \in
  \mathbb{F}_{q}\} \ .$$ 

  Similarly, we see in $\pi_{1}$ the Veronesean
  $$\mathcal{V}_{1} := [0,x_{1}^{2},0,x_{0}^{2},0,x_{0}x_{1},x_{1}x_{2},0,x_{2}^{2},x_{0}x_{2}]$$
  and the extra plane
  $$E_{1}:=\{[0,e_{0},0,0,0,e_{1},e_{2},0,0,0] \mid e_{0},e_{1},e_{2} \in
  \mathbb{F}_{q}\},$$ and in $\pi_{2}$, we have the Veronesean
  $$\mathcal{V}_{2} := [0,0,x_{2}^{2},0,x_{0}^{2},0,x_{1}^{2},x_{0}x_{2},x_{1}x_{2},x_{0}x_{1}]$$
  and the extra plane
  $$E_{2}:=\{[0,0,e_{0},0,0,0,0,e_{1},e_{2},0] \mid e_{0},e_{1},e_{2} \in
  \mathbb{F}_{q}\} \ .$$

Generalised dual arcs can be used to construct message authentication codes \cite{KLJ}. Below we give another application, namely secret sharing schemes.

\section{Secret sharing}\label{ss}

Now we will investigate applications of generalised dual arcs in
secret sharing schemes. For an overview of secret sharing and the
links with geometry we refer to \cite{Jackson:2004}. A recent overview
of different adversary models in secret sharing can be found in
\cite{Martin:2008}. 

We only consider a particular class of secret sharing schemes here, which is defined below.

\begin{Definition}
  In a \emph{$k$-out-of-$n$ secret sharing scheme} a \emph{dealer}
  generates $n$ shares $s_{1},\ldots,s_{n}$ and a secret $s$. The shares
  are given to different \emph{participants}. Each $k$ participants
  can reconstruct the secret with their shares.

  Less than $k$ participants cannot reconstruct the share. By $p_{i}$
  we denote the probability that $i<k$ participants may guess that
  share correctly. The probabilities $p_{i}$ are called the \emph{attack
    probabilities}. If $p_{i+1}/p_{i}>1$ the system \emph{leaks}
  information about the share.
\end{Definition}

Actually, we don't apply generalised dual arcs directly. But the dual
of these structures, which we call \emph{generalised arcs}.

\begin{Definition}
  A \emph{generalised arc} $\mathcal{A}$ of order $d$ with dimensions $n_1 <
  n_2 < \cdots < n_{d+1}$ of $PG(n,q)$ is a set of $n_{1}$-dimensional subspaces of $PG(n,q)$
  such that:
  \begin{enumerate}
  \item each $j$ of these subspaces generate a subspace of dimension
    $n_j$, $1 \leq j \leq d+1$,
  \item each $d+2$ of these subspaces span $PG(n,q)$.
  \end{enumerate}

  We call $(n,n_{1},\ldots,n_{d+1})$ the \emph{parameters} of the arc.

  If in addition the common intersection of all $n_{j+1}$-dimensional
  subspaces spanned by $j+1$ elements of the arc containing a given
  $n_{j}$-dimensional subspace $\pi$ spanned by $j$ elements of the arc
  is $\pi$, we call the arc \emph{regular}.
\end{Definition}

\begin{Theorem} \label{T:dual} The dual of an arc with parameters
  $(n,n_{1},\ldots,n_{d+1})$ is a dual arc with parameters
  $(n,n-1-n_{1},\ldots,n-1-n_{d+1})$ and vice versa.

Furthermore, the dual arc is regular if and only if the arc is
  regular.
\end{Theorem}
\begin{Proof}
  Dualising in $PG(n,q)$ maps every $k$-dimensional subspace onto an
  $(n-1-k)$-dimensional subspace. Dualising exchanges the concepts "span" and
  "intersection". \qed
\end{Proof}

Dual to Construction~\ref{cons:1}, we have the following construction
of generalised arcs.

\begin{Construction}\label{cons:2}
  As in Construction \ref{cons:1}, let $PG(V)$ be an $n$-dimensional
  space with basis $e_i$ ($0 \leq i \leq n$).

  Let $PG(W)$ be a $\left(\binom{n+d+1}{d+1}-1\right)$-dimensional space
  with basis $e_{i_0,\ldots,i_d}$ ($0 \leq i_0 \leq i_1 \leq \cdots \leq i_d \leq n$).

  We define $\zeta:PG(V) \to PG(W)$ by
  $$ \zeta: [\sum_{i=0}^{n} x_ie_i] \mapsto 
  [\sum_{0 \leq i_0 \leq \cdots \leq i_d \leq n} x_{i_0}\cdot\ldots\cdot x_{i_d} e_{i_0,\ldots,i_d}] \ .$$

  With $b$ and $B$ respectively, we denote the standard scalar product
  of $V$ and $W$, i.e.,
  $$ b(\sum_{i=0}^{n} x_i e_i,\sum_{i=0}^{n} y_i e_i)=\sum_{i=0}^{n} x_i y_i,$$  
  and
  $$ B(\sum_{0 \leq i_0 \leq  \cdots \leq i_d \leq n} x_{i_0,\ldots,i_d}
  e_{i_0,\ldots,i_d},\sum_{0 \leq i_0 \leq  \cdots \leq i_d \leq n} y_{i_0,\ldots,i_d} e_{i_0,\ldots,i_d})
  = \sum_{0 \leq i_0 \leq  \cdots \leq i_d \leq n} x_{i_0,\ldots,i_d}y_{i_0,\ldots,i_d}.
  $$

  For each $x \in V$, we denote by $x^{\perp}$ the subspace of $V$
  perpendicular to $x$ with respect to $b$. So
  $$ x^{\perp} = \{ y\in V \mid b(x,y)=0\}.$$

  For each point $P=[x]$ of $PG(V)$, we define a subspace $A(P)$ of
  $PG(W)$ by
  \begin{equation}
    \label{eq:A(P)}
    A(P) = \left< \zeta(y) \mid y \in x^{\perp} \right>.
  \end{equation}
\end{Construction}

\begin{Theorem}
  The set $\mathcal{A}=\{A(P) \mid P \in PG(n,q)\}$, defined in
  Construction~\ref{cons:2}, is a generalised arc with parameters
  $n_{i}=\binom{n+d+1}{d+1}-\binom{n+d+1-i}{d+1-i}-1$, $i=1,\ldots,d+1$.

  The generalised dual arc described in Construction \ref{cons:1} is
  the dual of that arc.
\end{Theorem}
\begin{Proof}
  By Definition (check Equation (\ref{eq:D(P)})), we have
  $D(P)=A(P)^{\bot}$ with respect to the bilinear form $B$. Since $B$ is a
  non-degenerate form, this means that $D(P)$ is dual to $A(P)$. Thus we
  may apply Theorem~\ref{T:dual}, which together with Theorem
  \ref{T:dual-arc} shows that $\mathcal{A}$ is indeed an arc. \qed
\end{Proof}

Before we describe the construction of a secret sharing scheme in
general, we give two examples that use the dual arc with parameters
$(9,5,2,0)$ we have seen in Example~\ref{E:2}.

\begin{Example} \label{E:s1} The dual of the dual arc with parameters
  $(9,5,2,0)$ is an arc consisting of $q^{2}+q+1$ different 
  $3$-dimensional spaces in $PG(9,q)$, with the
  following properties:
  \begin{enumerate}
  \item Each two $3$-dimensional spaces generate a $6$-space.
  \item Each three $3$-dimensional spaces generate an $8$-space.
  \item Each four $3$-dimensional spaces generate $PG(9,q)$.
  \end{enumerate}

  Now take the space $PG(10,q)$. Select any hyperplane as the secret.
  In that hyperplane select the above configuration of $q^{2}+q+1$
  $3$-dimensional spaces as shares.

  If the attacker does not have a share, he has a probability of
  $\frac{q-1}{q^{11}-1}$ to guess the secret $9$-space.

  If the attacker knows only one share, he has to guess a 9-space through
  the known $3$-dimensional space, so he has a probability of
  $\frac{q-1}{q^{7}-1}$ to guess the secret.

  Similarly, an attacker that knows $2$ or $3$ shares has a
  probability of $\frac{q-1}{q^{4}-1}$ or
  $\frac{q-1}{q^{2}-1}=\frac{1}{q+1}$ to guess the share.

  Any $4$ shares reconstruct the secret.
\end{Example}

\begin{Example} \label{E:s2}
  As in the previous example, we select a
  hyperplane $\Pi$ in $PG(10,q)$ and an arc consisting of $q^{2}+q+1$
  $3$-dimensional spaces with the same properties as above. One of
  these $3$-dimensional spaces
  $\pi$ will be the secret. 
  The other $3$-dimensional spaces are the shares.

  Furthermore, we select a $4$-dimensional space $\Pi_4$ through $\pi$ not
  contained in $\Pi$ 
  and make it public. If an attacker wants to find the secret space,
  he has to reconstruct $\Pi$ and then the secret space is the
  intersection $\Pi\cap\Pi_4$. 
  A short calculation shows that an attacker who
  knows $i$ ($i\leq4$) shares has a probability of
  $\frac{q-1}{q^{5-i}-1}$ to guess the secret.

  Another way to vary the attack probabilities is the following.
  Recall that the $q^{2}+q+1$ different $5$-spaces of the dual arc
  are of the form 
  $D(P)$ where $P$ is a point of a $2$-dimensional space
  $PG(2,q)$. The $q+1$ different
  $5$-spaces that correspond to the $q+1$ points of a line of $PG(2,q)$ lie in a
  common $8$-space. In the dual setting, this means that the $q+1$ corresponding
  $3$-dimensional spaces intersect in a common point.

  So if we fix one such $3$-dimensional space $\pi$, it has $q+1$ different intersection
  points with the other $q^{2}+q$ $3$-spaces. Suppose $\pi$ is the image of the point $P_{0}=[1,0,0]$.
  Furthermore, let $P_{1}=[0,1,0]$ and $P_{2}=[0,0,1]$. Consider lines of the form $\langle P_0,aP_{1}+P_{2}\rangle$ and $\langle P_0,P_1\rangle$.  
  Then they define $q+1$ different intersection points forming the
  twisted cubic arc consisting of the points $P_{a}=[1,a,a^{2},a^{3}]$ ($a \in
  \mathbb{F}_{q}$) and $P_{\infty}=[0,0,0,1]$. Choose a plane in $\pi$ which
  contains no intersection point. This is possible, since 
  $\mathbb{F}_{q}[X]$ contains an irreducible polynomial of degree
  $3$.

  Now we select this plane as the secret. We select  a $3$-dimensional space $\Pi_3$ through this plane not
  contained in $\Pi$ and make this public. An attacker who knows $i$ ($i<4$) shares has attack
  probabilities $p_{0}=p_{1}=\frac{1}{q^{3}+q^{2}+q+1}$,
  $p_{2}=\frac{1}{q^{2}+q+1}$ and $p_{3}=\frac{1}{q+1}$ to guess the
  secret. Thus the new scheme leaks no information if only one share
  is known.

  By selecting the correct subspace of $\pi$, we can also construct schemes
  that have no information leak for $2$ or $3$ shares. Then we must
  select a line or a point inside $\pi$  as the secret and take a plane $\Pi_2$ or line $\Pi_1$ through
  the selected line or point not in $\Pi$, and make this public.
\end{Example}

Now we give two theorems which use generalised arcs to construct
secret sharing schemes.

\begin{Theorem}
  In $PG(n+1,q)$, select an $n$-dimensional subspace $\Pi$ as the  secret.  In
  $\Pi$, select a generalised arc $\mathcal{A}$ of order $k-2$ with $n$
  elements and parameters $(n,d_{1},\ldots,d_{k-1})$. The elements of
  $\mathcal{A}$ are the shares.

  This describes a $k$-out-of-$n$ secret sharing scheme with the
  attack probabilities
  $$p_{i} = \frac{q-1}{q^{n+1-d_{i}}-1}$$
  for $0 \leq i < k$ (formally, we set $d_{0}=-1$).
\end{Theorem}
\begin{Proof}
  Every $k$ shares span $\Pi$, since $\mathcal{A}$ is a generalised arc of
  order $k-2$.

  Less than $k$ participants can take their shares $\pi_{1},\ldots,\pi_{i}$ and
  compute the $d_{i}$-dimensional space $\left<\pi_{1},\ldots,\pi_{i}\right>$.
  They know that $\Pi$ must contain that space. But for every
  $n$-dimensional space $\Pi'$ containing $\left<\pi_{1},\ldots,\pi_{i}\right>$,
  there exists an arc which has $\pi_{1},\ldots,\pi_{i}$ as elements. Thus the
  best attack is to guess an $n$-dimensional space through
  $\left<\pi_{1},\ldots,\pi_{i}\right>$. The number of such spaces is
  $\frac{q^{n+1-d_{i}}-1}{q-1}$. \qed
\end{Proof}

\begin{Theorem}
  In $PG(n+1,q)$, select a $(d_{1}+1)$-dimensional subspace $\pi'$ and
  make it public. In $\pi'$, select a $d_{1}$-dimensional subspace $\pi$ as
  the secret. Choose any hyperplane $\Pi$ of $PG(n+1,q)$ that contains $\pi$
  but not $\pi'$. Let $\mathcal{A}$ be a generalised dual arc of $\Pi$ of
  order $k-2$ with $n+1$ elements and parameters
  $(n,d_{1},\ldots,d_{k-1})$. The subspace $\pi$ should be an element of $\mathcal{A}$.
  The $n$ elements of $\mathcal{A}$ different from $\pi$ are the shares.

  This describes a $k$-out-of-$n$ secret sharing scheme with the
  attack probabilities
  $$p_{i} = \frac{q-1}{q^{d_{i+1}-d_{i}+1}-1}$$
  for $0 \leq i < k-1$ (formally, we set $d_{0}=-1$ and $d_{k}=n$).
\end{Theorem}
\begin{Proof}
  Every $k$ shares span $\Pi$, since $\mathcal{A}$ is a generalised arc of
  order $k-2$. Thus $k$ participants can compute $\Pi \cap \pi'$ which is the
  secret $\pi$.

  Less than $k$ participants can take their shares $\pi_{1},\ldots,\pi_{i}$ and
  compute the $d_{i}$-dimensional space $\left<\pi_{1},\ldots,\pi_{i}\right>$.
  Since the secret $\pi$ is also an element of the arc $\mathcal{A}$, we
  find that $\left<\pi_{1},\ldots,\pi_{i},\pi\right>$ has dimension $d_{i+1}$.
  This means that
  $$ \dim(\left<\pi_{1},\ldots,\pi_{i}\right> \cap \pi) = d_{i}+d_{1}-d_{i+1} \ .$$

  Since by construction $\pi' \cap \Pi = \pi$, we also have
  $$ \dim(\left<\pi_{1},\ldots,\pi_{i}\right> \cap \pi') =  d_{i}+d_{1}-d_{i+1} \ .$$

  The $i$ participants know that $\pi$ is a $d_{1}$-dimensional subspace
  of $\pi'$ containing the $(d_{i}+d_{1}-d_{i+1})$-dimensional subspace
  $\left<\pi_{1},\ldots,\pi_{i}\right> \cap \pi'$. But for every $d_{1}$-dimensional
  subspace $\bar{\pi}$ through $\left<\pi_{1},\ldots,\pi_{i}\right> \cap \pi'$ in
  $\pi'$, there exists a generalised arc containing 
  $\pi_{1},\ldots,\pi_{i}$ and $\bar{\pi}$. So the $i$ participants have no
  further information and must guess a $d_{1}$-dimensional subspace of
  $\pi'$ through $\left<\pi_{1},\ldots,\pi_{i}\right> \cap \pi'$. The probability for
  guessing this correctly is 
  $$p_{i}=\frac{q-1}{q^{d_{i+1}-d_{i}+1}-1} \ . $$
  \qed
\end{Proof}


\section{Known results}\label{kr}

In 1947, Bose studied ovals in \cite{Bose}. In that paper, he proved that 
an oval in $PG(2,q)$ has at most $q+1$ points if $q$ is odd and at most $q+2$ points if $q$ is even.

Special cases of generalised dual arcs have a long history. A
generalised dual arc of order $0$ is just a (partial) spread of
$PG(n,q)$. The generalised dual arc of order $n-1$ in $PG(n,q)$ with
parameters $(n,n-1,\ldots,1,0)$ is just the dual of an ordinary arc of
points in $PG(n,q)$.

Generalised dual arcs of order $1$ with $n_{2}=0$ are known as
$n_{1}$-dimensional dual arcs. It is known that the dimension $n$ of
the ambient space $PG(n,q)$ of an $n_{1}$-dimensional dual arc
satisfies $2n_{1} \leq n \leq \frac{1}{2}n_{1}(n_{1}+3)$ (see \cite{Yoshiara}).

\begin{Definition}
A family $\mathcal{A}$ of $\frac{q^{l+1}-1}{q-1}+1$ $l$-dimensional subspaces of $PG(n,q)$ with $n\geq 2$ is called an $l$-dimensional dual hyperoval if it satisfies the following three axioms:
\begin{itemize}
\item Every two elements of $\mathcal{A}$ intersect in a point.
\item Every three elements of $\mathcal{A}$ have no point in their intersection.
\item All members of $\mathcal{A}$ span the whole space $PG(n,q)$.
\end{itemize}
\end{Definition}


The next theorem is the translation to Veronesean dual arcs of the well-known fact that ovals in
$PG(2,q)$, $q$ odd, are maximal, but ovals in $PG(2,q)$, $q$ even, can
be extended to hyperovals  (see also~\cite{Yoshiara}).
\begin{Theorem} \label{T:d1-ext}
  For $q$
  odd, the Veronesean dual arc is maximal while for $q$ even, the Veronesean dual arc can
  be extended by an 
  $n_{1}$-dimensional space to an $n_{1}$-dimensional dual hyperoval. The
  extension element is called the {\em nucleus}.  
\end{Theorem}
\begin{Proof}
  In every arc element $\Omega=D([x_{0},\ldots,x_{n_1}])$, there is only
  one point not covered by a second arc element. This point is
  \[\zeta([x_{0},\ldots,x_{n_1}])=(x_{0}^{2},\ldots,x_{n_1}^{2},2x_{0}x_{1},\ldots,2x_{n_1-1}x_{n_1}),\] where $\zeta$ is the Veronesean map.
  
  For odd $q$, these points $\zeta([x_{0},\ldots,x_{n_1}])$ span $PG(\frac{1}{2}n_{1}(n_{1}+3),q)$, i.e. the
  Veronesean dual arc is not extendable. For $q$ even, they form an
  $n_{1}$-dimensional space which extends the Veronesean dual arc. This
  space is called the nucleus. \qed 
\end{Proof}

In 1958, Tallini \cite{Tallini:1958} (see also \cite{MR1363259}) showed that every
$2$-dimensional dual arc of $q^{2}+q+1$ elements in $PG(5,q)$, $q$ odd,
must be isomorphic to the dual arc defined by
Construction~\ref{cons:1}. This result was generalised 
in \cite{ThasHVM} to the following characterization of the finite
quadric Veronesean $\mathcal{V}_n^{2^n}$.
\begin{Result}\label{JefHendrik}
Let $\mathcal{F}$ be a set of $\frac{q^{n+1}-1}{q-1}$
$n$-dimensional spaces in $PG(\frac{n(n+3)}{2},q)$ with 
the following properties:
\begin{enumerate}[(1)]
\item[(VS1)] Each two elements of $\mathcal{F}$ intersect in a point.
\item[(VS2)] Each three elements of $\mathcal{F}$ are skew.
\item[(VS3)] The elements of $\mathcal{F}$ span $PG(\frac{n(n+3)}{2},q)$.
  \item[(VS4)] Any proper subspace of $PG(\frac{n(n+3)}{2},q)$ that is spanned by
  a collection of elements of $\mathcal{F}$ is a
  subspace of dimension $\frac{i(2n-i+3)}{2}-1$, for some
  $i \in \{0,\ldots,n\}$.
\item[(VS5)] If $q$ is even, at least one space spanned by two
  elements of $\mathcal{F}$ contains more than two elements of $\mathcal{F}$. 
\end{enumerate}
Then either $\mathcal{F}$ 
is a Veronesean dual arc with respect to
a quadric Veronesean $\mathcal{V}_n^{2^n}$ or $q$ is even and there are two members
$\Omega_1,\Omega_2\in\mathcal{F}$ such that the $2n$-dimensional space $\langle
\Omega_1,\Omega_2 \rangle$ only contains 2 elements of $\mathcal{F}$ and there
is a unique subspace $\Omega$ of dimension $n$ such that $\{\Omega\}\cup\mathcal{F}$ is
a Veronesean dual arc with the nucleus space as constructed in
Theorem~\ref{T:d1-ext}.
In particular, if $n=2$, then the statement holds under the weaker
hypotheses of $\mathcal{F}$ satisfying $(VS1)$, $(VS2)$, $(VS3)$
and $(VS5)$.     
\end{Result}

For order $d=1$ and $q$ even, there are non-Veronesean dual arcs with the
property that every space spanned by two elements of $\mathcal{F}$
contains exactly 
these two elements of $\mathcal{F}$. For $n=2$, one can classify all
examples that do 
not satisfy $(VS5)$ by a result of \cite{Delfra}; the only
possibilities are for $q=2$ and $q=4$. This classification remains open
for $n\geq 3$, although an infinite class of examples is known, described
in \cite{ThasHVM}.

\section{The case $d=1$} \label{sec:d1}

We prove that for $\delta>0$, $\delta$ small, a dual arc with parameters
$(n_{0},n_{1},n_{2})$ of size  $\frac{q^{n+1}-1}{q-1}-\delta$ is not maximal. The proof
techniques are similar to the techniques used in \cite{MR1363259} to give an
algebraic characterisation of a dual arc of size
$\frac{q^{n+1}-1}{q-1}$. The main difference is that the deficiency
$\delta$ makes simple counting arguments impossible, so we have to use more
difficult structural arguments.

\begin{Theorem} \label{T:d1}
Assume that $\delta \leq \frac{q-7}{2}$ for $q$ odd and $\delta \leq \frac{q-8}{2}$
for $q$ even, and let
$\mathcal{F}$ be a set of $\frac{q^{n+1}-1}{q-1}-\delta$ different
$n$-dimensional spaces in $PG(\frac{n(n+3)}{2},q)$ with 
the following properties:
\begin{enumerate}[(1)]
\item Each two elements of $\mathcal{F}$ intersect in a point.
\item Each three elements of $\mathcal{F}$ are skew.
\item The elements of $\mathcal{F}$ span $PG(\frac{n(n+3)}{2},q)$.
  \item Any proper subspace of $PG(\frac{n(n+3)}{2},q)$ that is spanned by
  a collection of elements of $\mathcal{F}$ is a
  subspace of dimension $\frac{i(2n-i+3)}{2}-1$, for some
  $i \in \{0,\ldots,n\}$.
\item If $q$ is even, at least one space spanned by two elements of
  $\mathcal{F}$ contains more than two elements of $\mathcal{F}$. 
\end{enumerate}
  Then $\mathcal{F}$ is
  extendable to a  regular generalised dual arc of size
  $\frac{q^{n+1}-1}{q-1}$. (In the case $q$ even, this dual arc of size $\frac{q^{n+1}-1}{q-1}$ is even extendable to a
  dual hyperoval.)
\end{Theorem}

The idea of the proof is in the same spirit as the proof of Result \ref{JefHendrik}, so the proofs of some results describing the general structure will look very similar as the ones used for that result.
The main work lies in the lemmata which actually deal with the deficiency itself, where we have to reconstruct the missing elements.

\begin{Definition}
A {\em contact point} is a point belonging to at most one element of $\mathcal{F}$.
\end{Definition}

Property (4) seems very technical. Our next lemma shows that for
large $q$, property (4) is no restriction. This motivates property (4).

\begin{Lemma} \label{l:line}
  Let $q \geq n$, then any configuration $\mathcal{F}$ which satisfies the
  properties (1)-(3) also satisfies property (4). 
\end{Lemma}
\begin{Proof}
  Assume that the claim of the
  lemma is wrong, i.e.\ there exists a sequence $\pi_{0},\ldots,\pi_{k}$ of
  elements in $\mathcal{F}$ with the property:

  \begin{itemize}
  \item $\Pi_{j} = \left<\pi_{0}, \ldots,\pi_{j}\right>$, for $j \leq k$,
  \item $\dim \Pi_{j} = \frac{(j+1)(2n-j+2)}{2}-1$, for $j<k$,
  \item $\frac{k(2n-k+3)}{2}-1 < \dim \Pi_{k} < \frac{(k+1)(2n-k+2)}{2}-1$.
  \end{itemize}

  By induction, we will construct a sequence $\pi_{k+1},\ldots,\pi_{n+1}$ of
  members of $\mathcal{F}$ with the properties:
  \begin{enumerate}[(I)]
  \item the subspace defined recursively by
    $\Pi_{i}=\left<\Pi_{i-1},\pi_{i}\right>$ has at least an $i$-dimensional
    subspace in common with $\pi_{i+1}$,
  \item the space $\pi_{i+1}$ is not contained in $\Pi_{i}$.
  \end{enumerate}
  For $i=n$, these two conditions yield a contradiction, because the elements of $\mathcal{F}$ have dimension $n$. This proves
  the lemma. 

  Now we construct $\pi_{j+1}$ from the sequence $\pi_{0},\ldots,\pi_{j}$.
  Note that $\dim \Pi_{j}$ is bounded by 
  \begin{multline*}
    \dim \Pi_{k}+(n-k)+\cdots+(n-(j-1)) \leq \\
    \frac{(k+1)(2n-k+2)}{2}-2 + \frac{(j-k)(2n-k-j+1)}{2} \leq \frac{n(n+3)}{2}-1
    \ .
  \end{multline*}
  Thus $\Pi_{j}$ is not the
  whole space. By property (3), we know that there exists a space $\bar{\pi}_{j+1}$ of $\mathcal{F}$
  not in $\Pi_{j}$. There are at least $q^{n}-1-\delta$ elements of
  $\mathcal{F}$ meeting $\bar{\pi}_{j+1}$ in a point outside of
  $\Pi_{j}$. Thus there are at least $q^{n}-\delta$ elements of $\mathcal{F}$
  not in $\Pi_{j}$. Since $\pi_{i+1}$ has at most an $(n-1)$-dimensional
  space in common with $\Pi_{i}$ ($i < k$), we conclude that at most
  $\frac{q^{n}-1}{q-1}$ elements of $\mathcal{F}$ intersect $\pi_{i+1}$
  in a point of $\Pi_{i}$. Thus for at most $j\frac{q^{n}-1}{q-1}$
  elements of $\mathcal{F}$, there exists an $i < j$ such that this
  element intersects $\pi_{i+1}$ in a point of $\Pi_{i}$. Because $k\leq n\leq q$, 
  $$q^{n}-\delta-k\frac{q^{n}-1}{q-1}>0,$$
  implying that there is an element $\pi_{j+1}$ of $\mathcal{F}$
  with the property that $\pi_{j+1}$ is not in $\Pi_{j}$ and $\pi_{j+1} \cap \pi_{i+1}
  \notin \Pi_{i}$. Especially, we have 
  $\dim \left< \pi_{j+1} \cap \pi_{i+1} \mid -1 \leq i<j \right>=j$, i.e. $\pi_{j+1} \cap
  \Pi_{j}$ is at least a $j$-dimensional space. 
  
  Thus, by induction, we have found the 
  members of $\mathcal{F}$ with the properties (I) and (II), which
  proves the lemma. \qed
\end{Proof}

Property (4) allows us to compute the dimensions of many objects
related to $\mathcal{F}$. An important special case is the following result.\\

{\bf Remark.} Let $\Pi$ be a $2n$-dimensional space spanned by two
  elements of $\mathcal{F}$. Then an element of $\mathcal{F}$
  either lies inside 
  $\Pi$ or intersects $\Pi$ in a line. \\  

The next lemma gives us an upper bound on the number of elements of
$\mathcal{F}$ contained in a space having one of the dimensions
mentioned in property (4).  

\begin{Lemma}\label{l:general}
Every $\left(\frac{i(2n-i+3)}{2}-1\right)$-dimensional space contains
at most $\frac{q^{i}-1}{q-1}$ elements of $\mathcal{F}$. 
\end{Lemma}
\begin{Proof}
Let $\Pi$ be an $\left(\frac{i(2n-i+3)}{2}-1\right)$-dimensional space
spanned by $i$ elements $\pi_{1}, \ldots, \pi_{i}$
of $\mathcal{F}$.

  An element of $\mathcal{F}$, not contained in $\Pi$, intersects $\Pi$ in
  an $(i-1)$-dimensional space $\Pi_i$ (this is part of property (4)).
  Each element of 
  $\mathcal{F}$, contained in $\Pi$, must share a 
  point with $\Pi_i$. Furthermore, no two elements of
  $\mathcal{F}$ in $\Pi$ intersect $\Pi_i$ in the same point, so $\Pi$
  contains at most $\frac{q^{i}-1}{q-1}$ elements of $\mathcal{F}$. \qed
\end{Proof}

To understand the goal of the next lemma, consider the dual arc obtained
by Construction \ref{cons:1}. In this example, every element of
$\mathcal{F}$ corresponds to a point of a projective space
$PG(n,q)$. The $2n$-dimensional spaces spanned by two elements of
$\mathcal{F}$ correspond to the lines of $PG(n,q)$. Thus if a dual arc
with $\frac{q^{n+1}-1}{q-1}-\delta$ elements is a subset of this example, then
the following is true: \\

\emph{Every $2n$-dimensional space spanned by two elements of
$\mathcal{F}$ contains at least $q+1-\delta$ elements of $\mathcal{F}$.}\\

Lemma \ref{l:q-d} is the first step in that direction.

\begin{Lemma} \label{l:q-d}
  Every $2n$-dimensional space contains $0$, $1$, $2$ or at least
  $q-\delta$ ($\delta\leq(q-7)/2$ for $q$ odd and $\delta \leq (q-8)/2$ for $q$ even) elements of $\mathcal{F}$. 

  If $q$ is odd, no $2n$-dimensional space contains exactly $2$ elements of
  $\mathcal{F}$.  
\end{Lemma}
\begin{Proof}
  Let $\Pi$ be a $2n$-dimensional space which contains $k$ elements of $\mathcal{F}$,
  where $2 \leq k < q-\delta$.
  
  Let $\pi'$ be any element of $\mathcal{F}$ not contained in $\Pi$.  
  This element $\pi'$ intersects $\Pi$ in a line $l'$ by the remark after Lemma \ref{l:line}. At
  least $q-\delta$ points of 
  $l'$ must be covered by a second element of $\mathcal{F}$. Since
  $q-\delta-k>0$, there must be a second element $\pi''$ of $\mathcal{F}$, not
  contained in $\Pi$, which intersects $l'$. Let $\pi'' \cap \Pi = l''$.
  
  The lines $l'$ and $l''$ span a plane $\pi$. Since every one of the $k$
  elements of $\mathcal{F}$ in $\Pi$ must intersect $\pi'$ and $\pi''$, these
  $k$ elements intersect $\pi'$ and $\pi''$ in a point on $l'$, respectively
  on $l''$, different from $l' \cap l''$. Hence, they intersect $\pi$ in lines.
  
  Assume that $\pi'''$ is another element of $\mathcal{F}$, not contained
  in $\Pi$, that intersects $\Pi$ in $l'''$. We prove that if $l'''$ has a
  point in common with $l'$, then it has also a point in common with
  $l''$.
  
  Suppose that $l'''$ intersects $l'$.  If $l'''$ does not intersect
  $l''$, then every element of $\mathcal{F}$ contained in $\Pi$ must share
  a line with the plane spanned by $l'$ and $l''$, and has a point in
  common with $l'''$. Thus these elements share a plane with the $3$-dimensional
  space spanned by $l'$, $l''$ and $l'''$. Especially, two of these elements intersect each other in a line, a contradiction.
  
  This proves that the elements of $\mathcal{F}$, not contained in $\Pi$,
  can be partitioned into {\em groups}. The elements from one group
  intersect each other in $\Pi$, and elements from different groups
  intersect each other outside of $\Pi$. Each group defines a plane inside
  $\Pi$ and the $k$ elements of $\mathcal{F}$ contained in $\Pi$ must
  intersect such a plane in lines.
  
  Let $\pi_1$ and $\pi_2$ be two planes inside $\Pi$ defined by such groups. We distinguish several cases for the intersection $\pi_1\cap\pi_2$.

  (1) The planes $\pi_{1}$ and $\pi_{2}$ cannot be skew to each
  other. Otherwise, they would span a $5$-dimensional space $\Omega$. Now
  every element of $\mathcal{F}$ in $\Pi$ shares a line with $\pi_{1}$ and
  $\pi_{2}$, so shares at least a $3$-dimensional space with $\Omega$, but
  then the elements of $\mathcal{F}$ in $\Pi$ intersect each other in at
  least a line, which is false. 

  (2) If $\pi_1$ and $\pi_2$ intersect in a line, then at most one element of
  $\mathcal{F}$ contained in $\Pi$ contains the line $\pi_1\cap \pi_2$. So at
  least $k-1$ elements of $\mathcal{F}$ contained in $\Pi$ must share a 
  plane with the
  $3$-dimensional space spanned by $\pi_1$ and $\pi_2$. Thus each two of
  these elements must share a line, a contradiction for $k>2$. We now
  eliminate the case $k=2$, where one of the two elements of
  $\mathcal{F}$ in $\Pi$, for instance $\pi$, passes through the line
  $\ell=\pi_1\cap \pi_2$.
  
  For $k=2$, all groups have size at least $q-\delta-1$. For, consider a
  first element $\pi'$ of $\mathcal{F}$ not in $\Pi$, then consider the line
  $\ell'=\pi'\cap \Pi$. This line has at most $\delta+1$ contact points, so it is
  intersected in a point by at least $q-2-\delta$ elements of
  $\mathcal{F}$, not lying in 
  $\Pi$. This shows that a group of elements of
  $\mathcal{F}$, not lying in $\Pi$, has at least size  $q-\delta-1$.
  
  But now consider the line $\ell=\pi_{1} \cap \pi_{2}$, lying in an element $\pi$ of
  $\mathcal{F}$ in $\Pi$, and in the two planes $\pi_1$ and
  $\pi_2$ containing at least $q-\delta-1$ lines lying in elements of
  $\mathcal{F}$, not contained in $\Pi$. Since no point of $\ell$ lies in three
  elements of $\mathcal{F}$, and every point of $\ell$ already lies in the
  element $\pi$ of $\mathcal{F}$, we must have $q+1\geq 2(q-\delta-1)+1$, where the
  $+1$ arises from the second element of $\mathcal{F}$ in $\Pi$. This
  implies $q\leq 2\delta+2$, a contradiction.

  (3) Thus $\pi_1$ and $\pi_2$ intersect in a point $Q$. But then the only
  possibility for an element of $\mathcal{F}$ contained in $\Pi$ to
  intersect $\pi_1$ and $\pi_2$ in lines is that $Q$ is a point of that
  element. Thus all elements of $\mathcal{F}$ contained in $\Pi$ contain
  $Q$. Since every three elements of $\mathcal{F}$ are skew, this means
  that $k=2$. 

  Assume now that we are in the case $k=2$ and $q$ is odd. Since there are
  $\frac{q^{n+1}-1}{q-1}-2-\delta$ elements of $\mathcal{F}$ not 
  contained in $\Pi$, and since for odd $q$ a dual arc of lines in
  $PG(2,q)$ contains at most $q+1$
  elements, each group can contain at most $q-1$ elements, so there are
  at least  
  $$\frac{1}{q-1}\left(\frac{q^{n+1}-1}{q-1}-2-\delta\right)>\frac{q^{n}-1}{q-1}$$ 
  different groups.
  
  Each group defines a plane through $Q$ which intersects an element of
  $\mathcal{F}$ contained in $\Pi$ in a line. Since an $n$-dimensional
  space only contains 
  $\frac{q^{n}-1}{q-1}$ different lines through $Q$, there must
  exist two groups which 
  define planes $\pi_1$ and $\pi_2$ intersecting in a line. But this is 
  impossible as we already proved. 
  
  So the case $k=2$ is only possible for $q$ even. \qed
\end{Proof}

Even if we could not exclude the case $k=2$ for $q$ even, we have
proven in step (3) the following characterisation:
\begin{Corollary} \label{C:q-d}
  Let $q$ be even and let $\left<\pi,\pi'\right>$ be a $2n$-space that contains only $\pi$
  and $\pi'$ as elements of $\mathcal{F}$. Then the elements of
  $\mathcal{F} \backslash \{\pi,\pi'\}$ intersect $\left<\pi,\pi'\right>$ in groups of
  pairwise intersecting lines. Furthermore, there can be at most
  $\frac{q^{n}-1}{q-1}$ such groups. 
\end{Corollary}

We call a $2n$-dimensional space {\em big} if it contains at least $q-\delta$
elements of $\mathcal{F}$. The next lemma associates with each big
$2n$-dimensional space $\Pi$ a plane $\bar{\pi}$ which will be very
important in the remaining part of this section.

\begin{Lemma} \label{l:bar-p}
  Let $\Pi$ be a $2n$-dimensional space containing $q+1-\delta_{i} \geq q-\delta$ elements of
  $\mathcal{F}$. Then $\Pi$ contains a plane $\bar{\pi}$ which intersects
  the $q+1-\delta_{i}$ elements of $\mathcal{F}$ in $\Pi$ in lines. The elements
  of $\mathcal{F}$, not in $\Pi$, intersect $\Pi$ in a line. These lines
  either lie in $\bar{\pi}$, or they are skew to $\bar{\pi}$ and then contain
  $\delta_{i}$ contact points. Moreover, those latter lines skew to $\bar{\pi}$ which are the intersection of 
   $\Pi$ with an element of $\mathcal{F}$ not lying in $\Pi$ are pairwise disjoint.
\end{Lemma}
\begin{Proof}
  Assume that two elements $\tilde{\pi}_{1}$ and $\tilde{\pi}_{2}$ of
  $\mathcal{F}$, not in 
  $\Pi$, intersect $\Pi$ in two intersecting lines $\ell_{1}$ and $\ell_{2}$. Let
  $\bar{\pi}$ be the plane spanned by $\ell_{1}$ and $\ell_{2}$.

  We are not in the case which is assumed in the beginning of the
  proof of Lemma \ref{l:q-d}. 
  However, the same kind of arguments as the ones used  in the
  proof of Lemma \ref{l:q-d} show that 
  \begin{enumerate}
  \item Every line in $\Pi$ 
    that intersects $\bar{\pi}$ and that comes from an element of
    $\mathcal{F}$ not in $\Pi$ must lie in $\bar{\pi}$.
  \item Every element of $\mathcal{F}$ in $\Pi$ must intersect $\bar{\pi}$
    in a line.
  \item The lines in $\Pi$ that come from an element of $\mathcal{F}$ not in $\Pi$
    and that do not lie in $\bar{\pi}$ must be pairwise disjoint.
  \end{enumerate}
  Property 3 is proven in the following way. Otherwise we have two
  planes $\bar{\pi}_{1}$ and $\bar{\pi}_{2}$ 
  corresponding with two different groups of lines as in the proof of
  Lemma \ref{l:q-d}. We have shown in 
  the proof of Lemma \ref{l:q-d} that $\bar{\pi}_{1}$ and $\bar{\pi}_{2}$ must
  intersect in a point $Q$ which lies on every element of $\mathcal{F}$
  in $\Pi$. But this implies that $\Pi$ has only $2$ elements of
  $\mathcal{F}$ which is not the case.

  So, from now on, we may assume that all the elements of $\mathcal{F}$, not in
  $\Pi$, intersect $\Pi$ in pairwise disjoint lines. Now we construct the
  plane $\bar{\pi}$.

  Let $\pi_{1}$, $\pi_{2}$ and $\pi_{3}$ be three elements of $\mathcal{F}$ in
  $\Pi$. Let $Q_{12}=\pi_{1} \cap \pi_{2}$, $Q_{13}=\pi_{1} \cap \pi_{3}$
  and $Q_{23}=\pi_{2} \cap \pi_{3}$.

  The points $Q_{12}$, $Q_{13}$, $Q_{23}$ generate a plane $\bar{\pi}$,
  since otherwise, 
  $\pi_1$, $\pi_2$, $\pi_3$ share a line. Assume that an element of
  $\mathcal{F}$, not in $\Pi$, intersects $\Pi$ in a line $\ell$ that meets
  $\bar{\pi}$. We claim that $\ell$ must lie in $\bar{\pi}$. Suppose
  the contrary. Without loss of generality, we may assume that $\ell \cap
  \bar{\pi} \notin \pi_{1} \cup \pi_{2}$. But then $\pi_{1}$ and $\pi_{2}$ share a plane
  with the
  $3$-dimensional space $\left<\bar{\pi},\ell\right>$, i.e. they share a
  line, a contradiction.

  At most one line in $\Pi$ that comes from an element of $\mathcal{F}$ not in
  $\Pi$ lies in $\bar{\pi}$, since these lines are pairwise disjoint. Since
  every element of $\mathcal{F}$ has only $\delta+1$ contact 
  points, this proves that at least $q-\delta-1$ points of $Q_{12}Q_{13}$
  lie in an element of $\mathcal{F}$ in $\Pi$, different from $\pi_{1}$.

  Assume that there exists an element $\pi$ of $\mathcal{F}$ in $\Pi$ which
  intersects $\pi_{1}$ in a point $Q$ not on $Q_{12}Q_{13}$. The above
  arguments show that $Q_{12}Q_{13}$, $QQ_{12}$ and $QQ_{13}$ must
  contain at least $3(q-\delta-1)-3>q+1$ points in $\pi_{1}$ which lie on two elements of
  $\mathcal{F}$ inside $\Pi$, a contradiction with Lemma \ref{l:general}.

  Thus every element $\pi$ of $\mathcal{F}$ in $\Pi$ meets $Q_{12}Q_{13}$,
  $Q_{12}Q_{23}$ and $Q_{13}Q_{23}$, i.e. it has a line in common with
  $\bar{\pi}$. \qed
\end{Proof}

The next series of lemmas deal with the case $q$ even and $k=2$. Let us
again have a look at the example that comes from Construction
\ref{cons:1}. In this example, every $2n$-dimensional space containing at least one element of
$\mathcal{F}$ contains
either $1$ or $q+1$ elements of $\mathcal{F}$. If $q$ is even, we can
extend the dual arc of size $\frac{q^{n+1}-1}{q-1}$ 
by one element $\pi$. This element $\pi$ has the special property
that for all other elements $\pi' \in \mathcal{F}$, the $2n$-space 
$\left<\pi,\pi'\right>$ contains no other element of $\mathcal{F}$, see
\cite{ThasHVM}.  We call this element the {\em nucleus} of $\mathcal{F}$.

We will prove in Lemma \ref{l:ext} that this property holds for every regular
generalised dual arc for $q$ even.

\begin{Lemma} \label{l:k2}
  Let $q$ be even and let $\pi,\pi' \in \mathcal{F}$ be such that the
  $2n$-dimensional space 
  $\left<\pi,\pi'\right>$ contains no other element of $\mathcal{F}$. Let
  $Q=\pi \cap \pi'$.

  Let $\Pi$ be a big $2n$-dimensional space containing $\pi$ and let $\bar{\pi}$ be the
  plane inside $\Pi$ described by Lemma \ref{l:bar-p}. Then $Q \in \bar{\pi}$.
\end{Lemma}
\begin{Proof}
  Let $\Pi = \left<\pi,\pi''\right>$, $\pi'' \in \mathcal{F} \backslash \{\pi,\pi'\}$. Let $\bar{\pi}'=\left<Q=\pi \cap \pi',\pi'' \cap \pi, \pi'' \cap
    \pi'\right>$. As we have already seen in Corollary~\ref{C:q-d}, this
  gives us a group of intersecting lines in this plane. But Lemma
  \ref{l:bar-p} states that the only plane in $\Pi$ which contains a group
  of intersecting lines is $\bar{\pi}$, i.e. $\bar{\pi}=\bar{\pi}'$. \qed.
\end{Proof}

\begin{Lemma} \label{l:k22}
  Let $q$ be even. For each $\pi \in \mathcal{F}$ either all $2n$-dimensional spaces
  $\left<\pi,\pi'\right>$ with $\pi \neq \pi' \in \mathcal{F}$ contain exactly two
  elements of $\mathcal{F}$, or there exists at most one element
  $\pi \neq \pi' \in \mathcal{F}$ such that $\left<\pi,\pi'\right>$ contains exactly two
  elements of $\mathcal{F}$.
\end{Lemma}
\begin{Proof}
  Assume that $\pi$ lies in a big $2n$-dimensional space $\Pi$, and let $\bar{\pi}$ be
  the plane described by Lemma \ref{l:bar-p} and let $\ell$ be the line
  $\bar{\pi} \cap \pi$. By Lemma \ref{l:k2}, we know that an element $\pi'$ of
  $\mathcal{F}$ for which $\left<\pi,\pi'\right>$ contains no other
  element of $\mathcal{F}$ must intersect $\pi$ in a point of $\ell$.

  Since $\ell$ has only $q+1$ points and
  $|\mathcal{F}|=\frac{q^{n+1}-1}{q-1}-\delta$, this means that $\pi$ must lie
  in more than one big 
  $2n$-space $\Pi'$. But then we have a second line $\ell' = \bar{\pi}' \cap \pi$ and
  every element $\pi'$ of
  $\mathcal{F}$ for which $\left<\pi,\pi'\right>$ contains no other
  element of $\mathcal{F}$ must intersect $\pi$ in a point of $\ell \cap
  \ell'$. ($\ell$ and $\ell'$ are different, since $\ell$ must meet the $q-\delta$ elements
  of $\mathcal{F}$ in $\Pi$, $\ell'$ must meet the $q-\delta$ elements
  of $\mathcal{F}$ in $\Pi'$ and $2q-2\delta-2>q+1$, see also step (2) of
  Lemma~\ref{l:q-d}.)
  This proves the lemma. \qed
\end{Proof}

\begin{Lemma} \label{l:ext}
  Let $q$ be even and assume that there exists a $2n$-dimensional space $\Pi$ which contains exactly two elements of $\mathcal{F}$. Then there exists at most one element $\pi \in \mathcal{F}$ such that
  for every $\pi \neq \pi' \in \mathcal{F}$, the $2n$-space $\left<\pi,\pi'\right>$
  contains exactly two elements of $\mathcal{F}$.
\end{Lemma}
\begin{Proof}
  Let $\Pi = \left<\pi,\pi'\right>$. Assume that both elements $\pi$ and $\pi'$ lie in a big $2n$-dimensional space. Then all other elements of $\mathcal{F}$ generate with
  $\pi$ and $\pi'$, respectively, a big $2n$-dimensional space (Lemma
  \ref{l:k22}). Let $\pi''$ be such 
  an element and $\Pi_{0}=\left<\pi,\pi''\right>$ with the special plane
  $\bar{\pi}_{0}$ and $\Pi_{1}=\left<\pi',\pi''\right>$ with the special plane
  $\bar{\pi}_{1}$. By the proof of Lemma \ref{l:k2}, we know that
  $\bar{\pi}_{0}=\left<\pi\cap \pi', \pi\cap \pi'', \pi'\cap \pi''\right>=\bar{\pi}_{1}$.

  But this is a contradiction since this plane cannot contain
  $2(q-\delta)-1>q+2$ different lines coming from elements of $\mathcal{F}$ in
  $\Pi_{0}$ and $\Pi_{1}$. Thus either $\pi$ or $\pi'$ does not lie in 
  big $2n$-dimensional spaces. They cannot both lie only in
  $2n$-spaces which contain $2$ elements of $\mathcal{F}$ or else by
  condition (5) of Theorem~\ref{T:d1} which we assume to be valid for $\mathcal{F}$, we find a $\pi'' \in
  \mathcal{F} \backslash \{\pi,\pi'\}$ 
  lying in at least one big $2n$-space and in
  two $2n$-spaces with only two 
  elements of $\mathcal{F}$, a contradiction with Lemma \ref{l:k22}. \qed
\end{Proof}

If $q$ is even and the special element $\pi$ from Lemma \ref{l:ext} exists, we
simply remove it from $\mathcal{F}$. This increases the deficiency by
$1$. \\

{\bf Remark.} Thus from now on, we assume that a $2n$-space cannot
contain $2$ elements of $\mathcal{F}$ and that $\delta \leq (q-6)/2$ when $q$ is even and $\delta \leq (q-7)/2$ when $q$ is odd.\\

Our next goal is a stronger version of Lemma \ref{l:bar-p} which states
that an element of $\mathcal{F}$, not in a big $2n$-space $\Pi$, must be skew to the
plane $\bar{\pi}$. We will reach this goal with Lemma \ref{l:bar-p2}.

\begin{Lemma} \label{l:dim3}
  Let $\Pi_{1}$, $\Pi_{2}$ and $\Pi_{3}$ be distinct $2n$-dimensional spaces containing at least
  $q-\delta$ elements of $\mathcal{F}$. 
  
  Then $\dim(\Pi_{1} \cap \Pi_{2} \cap \Pi_{3}) \leq n$.
\end{Lemma}
\begin{Proof}
  By property (4), we know that $\dim(\Pi_{1} \cap \Pi_{2})\leq n+1$.

  Assume that $\Pi_{1} \cap \Pi_{2} \cap \Pi_{3}$ is an $(n+1)$-dimensional space
  $\Pi$. Since two elements of $\mathcal{F}$ span a $2n$-dimensional
  space, the space $\Pi$ contains 
  at most one element of $\mathcal{F}$ and the other elements of
  $\mathcal{F}$ in $\Pi_{i}$ intersect $\Pi$ in a line.

  Let $\ell$ be such a line in $\Pi$ that comes from an element of $\mathcal{F}$
  in $\Pi_{1}$. The 
  elements in $\Pi_{2}$ and $\Pi_{3}$ intersect $\ell$ in a point. Since
  $2(q-\delta-1) > q+1$, some point of $\ell$ lies on an element of
  $\mathcal{F}$ in $\Pi_{1}$, 
  $\Pi_{2}$ and $\Pi_{3}$. A contradiction since each point lies on at most
  $2$ elements of $\mathcal{F}$. \qed
\end{Proof}

In the case of Construction \ref{cons:1}, we know that the big
$2n$-dimensional spaces correspond to the lines of an $n$-dimensional
projective space $PG(n,q)$. Thus in that case, 
every element of $\mathcal{F}$ lies in exactly $\frac{q^{n}-1}{q-1}$
big $2n$-spaces. Now we can prove this for a regular generalised dual arc.

\begin{Lemma} \label{l:2n}
   Let $\pi \in \mathcal{F}$. 
   Consider all $2n$-dimensional spaces through 
   $\pi$ containing at least $q-\delta$ elements of $\mathcal{F}$.
   Then the planes $\bar{\pi}$ of these $2n$-spaces intersect $\pi$ in
   different lines through a common point.

   Moreover, there are exactly $\frac{q^{n}-1}{q-1}$ different big $2n$-spaces
   through $\pi$.
\end{Lemma}
\begin{Proof}
  Let $\Pi$ and $\Pi'$ be two different $2n$-spaces through $\pi$, and let
  $\bar{\pi}$ and $\bar{\pi}'$ be the corresponding planes defined by Lemma \ref{l:bar-p}. By Lemma
  \ref{l:bar-p}, we know that $\pi \cap \bar{\pi}$ and $\pi \cap \bar{\pi}'$ are
  lines. These lines must be different since otherwise $\pi \cap \bar{\pi}=\pi \cap
  \bar{\pi}'$ would 
  contain at least $2(q-\delta-1)>q+1$ points lying on $\pi$ and on  
  another element of $\mathcal{F}$.

  By the proof of Lemma \ref{l:bar-p}, we know that at most $\delta+2$ elements
  of $\mathcal{F}$ not in $\Pi'$ intersect $\Pi'$ in lines contained in $\bar{\pi}'$. The other elements
  intersect $\Pi'$ in pairwise skew lines. Thus $\Pi$ contains at least
  $q-2\delta-3 \geq 3$ elements of $\mathcal{F}$ that intersect $\Pi'$ in
  pairwise skew lines; we call this set of lines
  $\mathcal{L}_{1}$. By symmetry, we know that $\Pi'$ contains
  at least
  $q-2\delta-3 \geq 3$ elements of $\mathcal{F}$ that intersect $\Pi$ in
  pairwise skew lines; we call this set of lines
  $\mathcal{L}_{2}$.

  Each line in $\mathcal{L}_{1}$ must intersect each line of $\mathcal{L}_{2}$,
  in the intersection point of the corresponding elements of
  $\mathcal{F}$. Thus $\mathcal{L}_{1}$ and $\mathcal{L}_{2}$ are the
  lines of two opposite reguli of a hyperbolic quadric $Q^{+}(3,q)$.  

  By Lemma \ref{l:bar-p}, we know that every element of $\mathcal{F}$
  in $\Pi$ has a line in common with $\bar{\pi}$. Thus the line $\pi \cap
  \bar{\pi}$ intersects all lines of $\mathcal{L}_{1}$, i.e. it lies in
  the regulus defined by $\mathcal{L}_{2}$. By symmetry, $\pi \cap \bar{\pi}'$
  lies in the regulus defined by $\mathcal{L}_{1}$. Thus $\pi \cap \bar{\pi}$
  and $\pi \cap \bar{\pi}'$ intersect. In addition we see that every element
  of $\Pi$ different from $\pi$ must lie in the regulus defined by
  $\mathcal{L}_{1}$, i.e. all elements of $\Pi$ intersect $\Pi'$ in
  pairwise skew lines not in $\bar{\pi}'$. Thus the first case in
  Lemma~\ref{l:bar-p} cannot occur. Especially the intersection point
  of $\pi \cap \bar{\pi}$ and $\pi \cap \bar{\pi}'$ must be a contact point, since
  it can lie only in elements of $\mathcal{F}$ that lie in the
  intersection $\Pi \cap \Pi'$. 

  This proves that either the lines of the form $\pi \cap \bar{\pi}$ share a
  common point or they lie in a common plane since they pairwise share a point. But the lines of the form
  $\pi \cap \bar{\pi}$ must additionally cover all non-contact points in $\pi$
  and intersect only in contact points. Thus 
  the lines of the form $\pi \cap \bar{\pi}$ share a 
  common contact point and there are at most $\frac{q^{n}-1}{q-1}$
  lines of the form $\pi \cap \bar{\pi}$. 
  That there are at least that many such
  lines follows from the fact that each big $2n$-space contains at most $q+1$
  elements of $\mathcal{F}$ and hence $\pi$ is contained in at least
  $(\frac{q^{n+1}-1}{q-1}-\delta-1)/q>\frac{q^{n}-1}{q-1}-1$ big $2n$-spaces.  
  \qed
\end{Proof}

{\bf Remark.} We note that in this proof, we encounter the strongest condition on $\delta$, namely $q-2\delta-3\geq 3$; equivalently, $\delta \leq (q-6)/2$ for $d=1$.\\

An important consequence of Lemma \ref{l:2n} is the following result.

\begin{Corollary}\label{cor:fulllines}
  Let $\Pi_{1}$, \ldots, $\Pi_{\frac{q^{n}-1}{q-1}}$ be the big
  $2n$-spaces containing a given element $\pi$ of $\mathcal{F}$. Let the space
  $\Pi_{i}$ contain $q+1-\delta_{i}$ elements of $\mathcal{F}$. Then
  $\sum_{i=1}^{\frac{q^{n}-1}{q-1}} \delta_{i} = \delta$.

  Especially, most big $2n$-dimensional spaces through $\pi$ contain $q+1$
  elements of $\mathcal{F}$. Moreover, each $2n$-space contains at least
  $q+1-\delta$ elements of $\mathcal{F}$. 
\end{Corollary}

\begin{Proof} We already know that every $2n$-space containing more than two elements of $\mathcal{F}$, contains $q+1-\delta_i\geq q-\delta$ elements of $\mathcal{F}$ (Lemma \ref{l:q-d}).

Since $\sum_{i=1}^{\frac{q^{n}-1}{q-1}} \delta_{i} = \delta$, necessarily $\delta_i\leq \delta$, so we can conclude that  every $2n$-space containing more than two elements of $\mathcal{F}$, contains $q+1-\delta_i\geq q+1-\delta$ elements of $\mathcal{F}$. \qed

\end{Proof}

The next lemma allows us to reduce the case of an
$(\frac{n(n+3)}{2},n,0)$-arc to the case of a $(5,2,0)$-arc.

\begin{Lemma} \label{l:plane}
  Let $\hat{\Pi}$ be a $(3n-1)$-space spanned by three elements of $\mathcal{F}$. Let
  $\hat{\mathcal{F}}$ be the set of elements of $\mathcal{F}$ in $\hat{\Pi}$. 

  For every $\pi$ in $\hat{\mathcal{F}}$, define
  $$ \hat{\pi} := \left<\pi \cap \pi' \mid \pi \neq \pi' \in \hat{\mathcal{F}}\right>
  .$$

  For every $\pi$ in $\hat{\mathcal{F}}$, the space $\hat{\pi}$ is a plane and these planes form
  a dual arc in $5$ dimensions. 
\end{Lemma}
\begin{Proof}
  For each element $\pi$ in $\hat{\mathcal{F}}$, we define a linear space
  $\mathcal{L}$ with
  the following properties:

  \begin{enumerate}[(i)]
  \item The points of the linear space are the points $\pi \cap \pi'$, with $\pi \neq \pi' \in
    \hat{\mathcal{F}}$.
  \item The lines of the linear space are the lines of $\pi$ through two
    points of the form $\pi \cap \pi'$ and $\pi \cap \pi''$ ($\pi'\neq\pi,\pi''\neq\pi \in
    \hat{\mathcal{F}}$). 

    If $\pi$ is not contained in the $2n$-dimensional space $\Pi$ spanned by $\pi'$ and $\pi''$, then it intersects $\Pi$ in a line containing  
    $\pi \cap \pi'$ and $\pi \cap \pi''$; in fact, this line contains at least $q-\delta$
    elements of the form $\pi \cap \pi'''$, with $\pi'''\neq\pi \in  \hat{\mathcal{F}}$
    (Lemma \ref{l:q-d}). 

    If $\pi$ is contained in the $2n$-space $\Pi$, then $\pi \cap \pi'$ and $\pi \cap
    \pi''$ lie on the intersection line of $\pi$ with the plane $\bar{\pi}$
    of $\Pi$ (Lemma \ref{l:bar-p}) which contains at least $q-\delta-1$
    intersection points of $\pi$ with other planes of $\hat{\mathcal{F}}$.
  \end{enumerate}
  
  The number of points in $\mathcal{L}$ is at least $3(q-\delta-1)-3$
  (Lemma \ref{l:q-d}) 
  and at most $q^{2}+q+1$ (Lemma \ref{l:general}).

  If $P_{0}$, $P_{1}$ and $P_{2}$ are three non-collinear points of
  the linear space, then
  $P_{0}P_{1}$ contains at least $q-\delta-1$ intersection points of two elements of $\mathcal{F}$ and thus there are at
  least $q-\delta-1$ lines through $P_{2}$ and therefore at least
  $(q-\delta-1)(q-\delta-2)+1$ intersection points in the plane
  $\left<P_{0},P_{1},P_{2}\right>$.

  By the same arguments, four points $P_{0}$, $P_{1}$, $P_{2}$ and $P_{3}$
  of the linear space $\mathcal{L}$
  that do not lie  in a plane would imply that the linear space
  $\mathcal{L}$ contains 
  at least $(q-\delta-2)[(q-\delta-1)(q-\delta-2)+1]+1$ points. But this is not
  possible since the number of points in $\mathcal{L}$ is bounded by
  $q^{2}+q+1$. 

  Thus $\hat{\pi} := \left<\pi \cap \pi' \mid \pi \neq \pi' \in \hat{\mathcal{F}}\right>$
  is a plane.

  It remains to be proven that the planes $\hat{\pi}$ span a
  $5$-space. Their span has at most dimension $5$ as we know from
  \cite{Yoshiara}. Assume that they only span a $4$-space. Then the three
  elements of $\mathcal{F}$ that span $\hat{\Pi}$ would have a plane in common
  with this $4$-dimensional space. This would imply that $\hat{\Pi}$ has at
  most dimension 
  $4+3(n-2) = 3n-2$, but this is false. 
  \qed
\end{Proof}
\begin{Corollary}\label{cor:2n-3n}
  Every big $2n$-space lies in exactly $\frac{q^{n-1}-1}{q-1}$
  different $(3n-1)$-spaces spanned by three elements of $\mathcal{F}$.
\end{Corollary}
\begin{Proof}
  Every big $2n$-space $\Pi$ through $\pi$ corresponds to a line $\bar{\pi} \cap \pi$
  and all these lines go through a common point $Q$. A $(3n-1)$-space defined by three elements of $\mathcal{F}$
  through $\pi$ corresponds to the lines in one plane through $Q$ inside $\pi$ by Lemma \ref{l:plane}. Thus
  the $(3n-1)$-spaces defined by three elements of $\mathcal{F}$ through $\Pi$ correspond to planes inside
  $\pi$ through $\bar{\pi} \cap \pi$. There are exactly $\frac{q^{n-1}-1}{q-1}$
  such planes. \qed
\end{Proof}

Now we are able to improve the result of Lemma \ref{l:bar-p}.

\begin{Lemma} \label{l:bar-p2}
  With the notations of Lemma \ref{l:bar-p}, the following result holds:
  No element of $\mathcal{F}$ not in $\Pi$ intersects $\bar{\pi}$.
\end{Lemma}
\begin{Proof}
  We know from the preceding lemma that $\bar{\pi}$ shares a line
  with every element $\pi$ of $\mathcal{F}$ in $\Pi$, passing through a fixed
  contact point of $\pi$. Assume that $\bar{\pi}$ 
  contains an extra line from an element $\pi'$ of $\mathcal{F}$ not contained
  in $\Pi$.  Let $\{R\}=\pi \cap \pi' \cap \bar{\pi}$.

  The elements $\pi$ and $\pi'$ define a big $2n$-dimensional space $\Pi'$,
  and $\Pi'$ contains 
  a plane $\bar{\pi}'$. The intersection $\bar{\pi}' \cap \pi$ is a line which
  contains $R$ and 
  the fixed contact point. Thus $\pi \cap \bar{\pi}' = \pi \cap \bar{\pi}$, a
  contradiction. 

  Thus $\bar{\pi}$ contains no line that comes from an element of
  $\mathcal{F}$ not in $\Pi$. Elements of $\mathcal{F}$ inside $\Pi$ 
  intersect $\bar{\pi}$ in a dual arc of $q+1-\delta_{i}$ lines. 
  \qed 
\end{Proof}

{\bf Remark.} If $\bar{\pi}$ contains lines of contact points, these lines
extend the dual arc of  $q+1-\delta_{i}$ lines induced by the elements of
$\mathcal{F}$ in $\Pi$. For $\delta_{i}=1$ and $q$ odd, we find
one line of contact points, and for $\delta_{i}=1$ and $q$ even, we find two
lines of contact points.\\

Now we are reaching our final goal to prove that $\mathcal{F}$ is not
maximal. As a first step, we prove that the planes $\bar{\pi}$ contain
lines of contact points.

\begin{Lemma} \label{l:intersect}
  Let $\Pi_1$ and $\Pi_2$ be two big $2n$-spaces with the property that
  $\left<\Pi_{1},\Pi_{2}\right>$ is a $(3n-1)$-dimensional space. 
  Assume that $\Pi_{1}$ and $\Pi_{2}$ share no
  element of $\mathcal{F}$. Let $\bar{\pi}_1$ and 
  $\bar{\pi}_2$ be the planes in $\Pi_1$ and $\Pi_2$ which exist by Lemmas
  \ref{l:bar-p} and \ref{l:bar-p2}. Then $\bar{\pi}_{1} \cap \Pi_{2}$ is a line of
  contact points.
\end{Lemma}

\begin{Proof}
  First of all, it is impossible that the plane $\bar{\pi}_1$ is
  contained in $\Pi_2$. For assume the contrary. We obtain a
  contradiction in the following way. Every element $\pi$ of $\mathcal{F}$
  in $\Pi_1$ intersects $\Pi_2$ in a line. If $\bar{\pi}_1$ lies completely
  in $\Pi_2$, then the intersection line $\ell=\Pi_2\cap\pi$ equals the line
  $\bar{\pi}_1\cap \pi$. This line contains at least $q-\delta_{1}$ points lying in
  two elements of $\mathcal{F}$ in $\Pi_1$. But the $q+1-\delta_{2}$ elements of
  $\mathcal{F}$ in $\Pi_2$ must intersect $\pi$ in a point. So at least
  $q+1-\delta_{2}$ points of $\ell$ still lie in an element of $\mathcal{F}$ in $\Pi_2$.
  Then there are points of $\ell$ lying in three elements of $\mathcal{F}$.
  This is false.

  Note that the plane $\bar{\pi}_{1}$ lies in the $5$-space $\hat{\Pi}_{1}
  \subseteq \Pi_{1}$
  spanned by the 
  planes $\hat{\pi}$ defined in Lemma \ref{l:plane}. 
  Then $\Pi_{2}$ cannot contain $\hat{\Pi}_{1}$, since otherwise every element of
  $\hat{\mathcal{F}}$ would intersect $\Pi_{2}$ at least in a plane,
  contradicting the remark after Lemma \ref{l:line}.
  Thus $\Pi_{2} \cap \hat{\Pi}_{1}$ is
  a $4$-dimensional space, spanned by two planes $\hat{\pi}$ and
  $\hat{\pi}'$ corresponding to elements $\pi$ and $\pi'$ of $\mathcal{F}$
  in $\Pi_{2}$. The plane $\bar{\pi}_{1}$ lies in the $5$-dimensional space
  $\hat{\Pi}_{1}$ and thus it intersects the $4$-dimensional space
  $\Pi_{2} \cap \hat{\Pi}_1$, 
  and therefore $\Pi_{2}$, in at least a line. 

  Consider again the intersection line $\ell=\Pi_2 \cap \pi$ of an element $\pi$ of
  $\mathcal{F}$ in $\Pi_1$ with $\Pi_{2}$. This line contains $q+1-\delta_{2}$
  points lying on an element  of $\mathcal{F}$ in $\Pi_2$ and $\delta_{2}$
  contact points (Lemma \ref{l:bar-p} and Lemma \ref{l:bar-p2}). So the points
  of $\ell$ do not lie in an 
  other element of $\mathcal{F}$ in $\Pi_1$.

  Now $\ell$ and $\pi \cap \bar{\pi}_1$ intersect in a point, since both lines
  lie in the plane $\hat{\pi}$ defined by Lemma \ref{l:plane}. This point must be
  a contact point, for else, it lies in a second plane of
  $\mathcal{F}$ in $\Pi_1$, but this was excluded in the preceding 
  paragraph.

  So $\pi$ shares a contact point with $\Pi_2$, which also lies on the
  intersection line of $\Pi_2$ with $\bar{\pi}_1$. 

  This proves that the line $\bar{\pi}_{1} \cap \Pi_{2}$ intersects the dual
  arc in $\bar{\pi}_{1}$, consisting of lines of the form $\bar{\pi}_{1} \cap
  \pi$, where $\pi$ is an element of $\mathcal{F}$ in $\Pi_{1}$, only in contact
  points, i.e. $\bar{\pi}_{1} \cap \Pi_{2}$ only contains points covered by
  at most one element of $\mathcal{F}$.

%
   
  This proves the theorem.
  \qed
\end{Proof}

\begin{Lemma} \label{l:regulus}
  Let $\Pi_1$ and $\Pi_2$ be two big $2n$-dimensional spaces with the property that
  $\left<\Pi_{1},\Pi_{2}\right>$ is a $(3n-1)$-space. 
  Assume that $\Pi_{1}$ and $\Pi_{2}$ share no
  element of $\mathcal{F}$ and let $q+1-\delta_{1}$ be the number of
  elements of $\mathcal{F}$ in $\Pi_{1}$ and let $q+1-\delta_{2}$ be the number
  of elements of $\mathcal{F}$ in $\Pi_{2}$.

  Let $\bar{\pi}_{1}$ and $\bar{\pi}_{2}$ be the planes in $\Pi_{1}$ and
  $\Pi_{2}$ which exist by 
  Lemma \ref{l:bar-p}.

  Then the lines $\mathcal{L}_{1}= \{ \bar{\pi}_{2} \cap \Pi_{1} \} \cup \{ \pi_{1} \cap
  \Pi_{2} \mid \pi_{1} \in \mathcal{F},\pi_{1} \subset \Pi_{1}\}$ and
  $\mathcal{L}_{2}= \{ \bar{\pi}_{1} \cap \Pi_{2} \} \cup \{ \pi_{2} \cap
  \Pi_{1} \mid \pi_{2} \in \mathcal{F},\pi_{2} \subset \Pi_{2}\}$ are lines of two
  opposite reguli of a hyperbolic quadric $Q^{+}(3,q)$.

  Especially this implies that $\delta_{1}>0$ and $\delta_{2}>0$ since a regulus
  has only $q+1$ lines, and that $\bar{\pi}_{2} \cap \Pi_{1}$ and
  $\bar{\pi}_{1} \cap \Pi_{2}$ are concurrent.
\end{Lemma}
\begin{Proof}
  By Lemma \ref{l:bar-p2}, we know that the elements of $\mathcal{F}$
  in $\Pi_{1}$ intersect $\Pi_{2}$ in pairwise skew lines. Thus
  $\mathcal{L}'_{1}= \{ \pi_{1} \cap \Pi_{2} \mid \pi_{1} \in \mathcal{F},\pi_{1} \subset \Pi_{1}\}$ and
  $\mathcal{L}'_{2}= \{ \pi_{2} \cap \Pi_{1} \mid \pi_{2} \in \mathcal{F},\pi_{2} \subset
  \Pi_{2}\}$ are sets of pairwise skew lines. Since $\pi_{1} \cap \pi_{2} \subset
  \Pi_{1} \cap \Pi_{2}$, every line of $\mathcal{L}'_{1}$ intersects every
  line of $\mathcal{L}'_{2}$.

  Since both sets contain more than $2$ lines, it follows that the lines
  of  $\mathcal{L}'_{1}$ and $\mathcal{L}'_{2}$ are lines of opposite
  reguli.

  Now consider the line $\bar{\pi}_{2} \cap \Pi_{1}$, which exists by Lemma
  \ref{l:intersect}. By Lemma \ref{l:bar-p2}, the plane $\bar{\pi}_{2}$
  is skew to all elements of $\mathcal{F}$ in $\Pi_{1}$. Thus
  $\bar{\pi}_{2} \cap \Pi_{1}$ is different from all lines in
  $\mathcal{L}'_{1}$. But every element $\pi_{2}$ of $\mathcal{F}$, contained in
  $\Pi_{2}$, has a line in common with $\bar{\pi}_{2}$. Thus $\bar{\pi}_{2} \cap
  \Pi_{1}$ intersects all lines of $\mathcal{L}'_{2}$. This proves that
  $\mathcal{L}_{1}= \{ \bar{\pi}_{2} \cap \Pi_{1} \} \cup \mathcal{L}_{1}'$ are
  the lines of a regulus. By symmetry, the same is true for
  $\mathcal{L}_{2}$. \qed
\end{Proof}

Recall that the final goal is to prove that $\mathcal{F}$ is given
by Construction \ref{cons:1}. Thus every element of $\mathcal{F}$
should correspond to a point of $PG(n,q)$. Since $\mathcal{F}$ has
only $\frac{q^{n+1}-1}{q-1}-\delta$ elements, $\delta$ points of $PG(n,q)$ are
not used in Construction \ref{cons:1}. The next lemma will identify
these \emph{holes}. 

Consider the linear space $\mathcal{L}$ with the elements of
$\mathcal{F}$ as points and the $2n$-spaces generated by two elements
of $\mathcal{F}$ as lines. This is a linear
space with $\frac{q^{n+1}-1}{q-1}-\delta$ points.

As planes of
$\mathcal{L}$, we define the $(3n-1)$-dimensional spaces generated by three elements of
$\mathcal{F}$.

\begin{Lemma} \label{l:hole}
  Every plane of $\mathcal{L}$ is a projective plane of order $q$ with
  possibly some holes.
\end{Lemma}
\begin{Proof}
  Let $P$ be a $(3n-1)$-dimensional space generated by three elements of
  $\mathcal{F}$. 

  Let $\Pi \subset P$ be a $2n$-dimensional space that contains $q+1$ elements of
  $\mathcal{F}$. This $2n$-space $\Pi$ exists since there are $q+1$ different
  big $2n$-dimensional spaces through an element $\pi$ of $\mathcal{F}$ in
  $P$ and at most $\delta$ of them contain less than $q+1$ elements of
  $\mathcal{F}$ (Corollary~\ref{cor:fulllines}). Let $\Pi'$ be an other
  big $2n$-dimensional space in $P$.  
  By Lemma \ref{l:regulus}, we know that $\Pi$ and $\Pi'$ share an element
  of $\mathcal{F}$.

  
  Let $\Pi_{1}$ and $\Pi_{2}$ be two big $2n$-dimensional spaces in
  $P$. Let $\pi$ be an element of $\mathcal{F}$ in $\Pi_{1}$, but not in
  $\Pi_{2}$. Since $\Pi_{2}$ contains at least $q+1-\delta$ elements of
  $\mathcal{F}$, there must be at least $q+1-\delta$ big $2n$-dimensional spaces in $P$
  through $\pi$.

  At most $\delta$ of the $\frac{q^{n}-1}{q-1}$ different big
  $2n$-dimensional spaces
  through $\pi$ contain less than $q+1$ elements of $\mathcal{F}$ (Corollary~\ref{cor:fulllines}).
  Each of the at least $q+1-\delta$ elements of $\Pi_{2}$ spans together with
  $\pi$ a big $2n$-dimensional space in $P$.
  Thus there are at least 
  $q+1-2\delta \geq 2$ big $2n$-spaces in $P$ through $\pi$ which contain exactly $q+1$
  elements of $\mathcal{F}$. We denote these $2n$-dimensional spaces
  by $\Pi_{1}$ and  $\Pi_{2}$. 

  By the same arguments we find an additional $2n$-dimensional space $\Pi_{3}$ which
  contains $q+1$ elements of $\mathcal{F}$, and which intersects $\Pi_{1}$ and
  $\Pi_{2}$ in different elements of $\mathcal{F}$.

  Thus that plane $P$ of $\mathcal{L}$ contains a triangle $\Pi_{1}$,
  $\Pi_{2}$, $\Pi_{3}$, and each side of the triangle contains $q+1$
  points of $\mathcal{L}$. Every other big $2n$-dimensional space in
  $P$ intersects $\Pi_{1}$, 
  $\Pi_{2}$ and $\Pi_{3}$ in elements of $\mathcal{F}$
  (Lemma~\ref{l:regulus}), thus a direct 
  counting argument shows us that $P$ contains
  $(q-1)^{2}+3(q-1)+3=q^{2}+q+1$ lines of $\mathcal{L}$, where
  $(q-1)^{2}$ is the number of lines intersecting the side of the
  triangle in different points, $3(q-1)$ is the number of lines
  through a vertex different from the sides and $3$ is the number of
  sides of the triangle. The number of
  elements of $\mathcal{F}$ in $P$ is at most $q^{2}+q+1$ (by Lemma
  \ref{l:general}) and 
  at least $q^{2}+q+1-\delta$ (by Corollary \ref{cor:fulllines}).

  Now consider any line $\ell$ of $P$ with $q+1-x$ points ($x \geq
  1$). Then $q(q+1-x)$ 
  lines intersect $\ell$ and thus there are $xq$ lines skew to $\ell$. Every
  point not on $\ell$ lies on $x$ lines that do 
  not intersect $\ell$. Thus there must exist a point not on $\ell$ that
  lies on a line $\ell'$ disjoint to $\ell$ with at least
  $\frac{[q^{2}+q+1-\delta-(q+1-x)]x}{qx}>q-1$ points. By 
  Lemma~\ref{l:regulus}, $\ell'$ has $q$ points. 

  As we have seen above, there are $q$ lines of $P$ skew to $\ell'$ and every
  point not on $\ell'$ lies on such a line. We may extend $\mathcal{L}$ by
  a point that lies on $\ell'$ and all lines skew to $\ell'$. Extending
  $\mathcal{L}$  stepwise by at most $\delta$ points, we obtain a
  $2-(q^{2}+q+1,q+1,1)$ design, i.e. a projective plane
  of order $q$.
  \qed
\end{Proof}

\begin{Lemma} \label{l:fine}
  Let $\mathcal{F}$ be a dual arc that satisfies the assumptions of
  Theorem~\ref{T:d1}. Let $\delta>0$, then $\mathcal{F}$ is not maximal.
\end{Lemma}
\begin{Proof}
  Since $\delta>0$, we find a big $2n$-space which contains less than $q+1$
  elements of $\mathcal{F}$ (Corollary~\ref{cor:fulllines}). Every
  $(3n-1)$-dimensional space spanned by three elements of $\mathcal{F}$
  through such a $2n$-space contains less 
  than $q^{2}+q+1$ elements of $\mathcal{F}$. Let $P$ be such a $(3n-1)$-space.

  Select a $2n$-space $\Pi$ in $P$ that contains exactly $q$ elements of
  $\mathcal{F}$. Such a space exists, because by Lemma \ref{l:hole},
  the linear space $\mathcal{L}$ is a projective plane with at most $\delta$ holes and
  such linear spaces contain lines with exactly $q$ points.

  Consider the $(3n-1)$-spaces through $\Pi$ generated by three elements
  of $\mathcal{F}$. By Lemma \ref{l:hole}, these $(3n-1)$-spaces define 
  projective planes with 
  holes. We will call a big $2n$-space $\Pi'$ parallel to $\Pi$ if it 
  "goes through" the unique hole of $\Pi$ in the corresponding
  projective plane defined by Lemma \ref{l:hole}.

  The $2n$-spaces parallel to $\Pi$ partition the set
  $\mathcal{F}$. By Corollary~\ref{cor:2n-3n}, we know
  that every big $2n$-space lies in $\frac{q^{n-1}-1}{q-1}$ different
  $(3n-1)$-spaces spanned by three elements of $\mathcal{F}$. Thus there
  are exactly
  $$ q\left(\frac{q^{n-1}-1}{q-1}\right)+1= \frac{q^{n}-1}{q-1}$$
  $2n$-spaces parallel to $\Pi$, including $\Pi$ itself.

  Consider two big $2n$-spaces $\Pi_{1}$ and $\Pi_{2}$ parallel to $\Pi$. If
  $\Pi_{2} \not\subset \left<\Pi,\Pi_{1}\right>$, then $\left<\Pi,\Pi_{1},\Pi_{2}\right>$ is
  a $(4n-3)$-dimensional space (Property (4)). Since $2n < \dim
  \left<\Pi_{1},\Pi_{2}\right> < \dim\left<\Pi,\Pi_{1},\Pi_{2}\right>=4n-3$,
  Property (4) implies that $\dim \left<\Pi_{1},\Pi_{2}\right>$ $ =3n-1$. Thus any
  two elements in the parallel class satisfy the conditions of Lemma
  \ref{l:intersect} and Lemma \ref{l:regulus}, i.e. they lie in a $(3n-1)$-space. 

  Let $q$ be odd. Choose any $2n$-space $\Pi'$ parallel to $\Pi$ which
  contains exactly $q$ elements of $\mathcal{F}$. By a direct counting
  argument we find that at least $\frac{q^{n}-1}{q-1}-(\delta-1)$ of the
  $\frac{q^{n}-1}{q-1}$ elements in the parallel class have this
  property. Then by Lemma~\ref{l:bar-p2}, the plane
  $\bar{\pi}'$ of $\Pi'$ contains exactly one line of contact points. 
  By Lemma \ref{l:intersect}, these lines must lie in the common intersection
  $\Omega$ of all $2n$-spaces parallel to $\Pi$. Thus $\Omega$ contains
  $\frac{q^{n}-1}{q-1}-(\delta-1)$ lines of contact points that share a common
  point $Q$ (Lemma \ref{l:regulus}). This proves that $\Omega$ is an
  $n$-dimensional space; it cannot be bigger by Lemma
  \ref{l:dim3}. Now look at any big $2n$-space $\Pi''$ parallel to
  $\Pi$ containing $q+1-\delta_{i}$ elements of $\mathcal{F}$. By Lemma
  \ref{l:intersect} and Lemma \ref{l:regulus}, the plane 
  $\bar{\pi}''$ must share a line through $Q$ with every other
  $2n$-space parallel to $\Pi$. This line must lie in $\Omega$ since
  otherwise $\bar{\pi}''$ would need different lines for each
  $2n$-space.  
  Thus $\Omega$ contains $\frac{q^{n}-1}{q-1}$ lines of contact points
  through $Q$, i.e., it only contains contact points and we can extend
  $\mathcal{F}$ by $\Omega$.

  For $q$ even, the situation is more complicated. We have always two
  lines of contact points and we must choose the correct one. Let $\Pi_1$ be
  a $2n$-space which contains exactly $q$ elements of
  $\mathcal{F}$. By  Lemma~\ref{l:intersect}, the plane $\bar{\pi}_{1}$ of
  $\Pi_{1}$ must share a line of contact points with each $2n$-space
  parallel to $\Pi_{1}$. By the pigeon hole principle there are at least
  $\frac{1}{2}[\frac{q^{n}-1}{q-1}-\delta]$ different $2n$-spaces
  parallel to $\Pi_{1}$, which contain $q$ elements of $\mathcal{F}$ and
  which intersect $\bar{\pi}_{1}$ in the same line $\ell_{1}$ of contact points.

  Let $\Pi_{2}$ and $\Pi_{3}$ be two such spaces. Choose $\Pi_{2}$ and $\Pi_{3}$ such
  that $\dim (\Pi_{1} \cap \Pi_{2} \cap \Pi_{3})=n$. For $n=2$, this is always the case since
  the intersection of three $4$-spaces in a $5$-space is at least a
  plane, and since Lemma \ref{l:dim3} states that $\dim (\Pi_{1} \cap
  \Pi_{2} \cap \Pi_{3})\leq 2$. For $n>2$, we can choose $\Pi_{2}$ and $\Pi_{3}$ such that 
  $\dim \left<\Pi_{1} , \Pi_{2} , \Pi_{3}\right>=4n-3$ and then we obtain $\dim (\Pi_{1} \cap
  \Pi_{2} \cap \Pi_{3})=n$ by the dimension formula.

  Let $\ell_{2}$ be the line $\bar{\pi}_{2} \cap \Pi_{1}$. Consider the
  hyperbolic quadric with the two reguli 
  $\mathcal{L}_{1}= \{ \bar{\pi}_{2} \cap \Pi_{1} \} \cup \{ \pi_{1} \cap
  \Pi_{2} \mid \pi_{1} \in \mathcal{F},\pi_{1} \subset \Pi_{1}\}$ and
  $\mathcal{L}_{2}= \{ \bar{\pi}_{1} \cap \Pi_{2} \} \cup \{ \pi_{2} \cap
  \Pi_{1} \mid \pi_{2} \in \mathcal{F},\pi_{2} \subset \Pi_{2}\}$ (see
  Lemma~\ref{l:regulus}).

  Then $\Pi_{3}$ contains the line $\ell_{1}=\bar{\pi}_1\cap \Pi_2$ of this hyperbolic quadric  since $\Pi_1$ shares the same line of $\bar{\pi}_1$ with $\Pi_2$ and $\Pi_3$.
  Hence, $\Pi_3$ must contain a second line of this hyperbolic quadric. We prove this as follows. We know that $\dim(\Pi_1\cap \Pi_2)=n+1$ and that $\dim (\Pi_{1} \cap
  \Pi_{2} \cap \Pi_{3})=n$. The hyperbolic quadric $\mathcal{L}_{1} \cup \mathcal{L}_{2}$ cannot lie in  $\Pi_{1} \cap
  \Pi_{2} \cap \Pi_{3}$, or else every space $\pi_1\in \mathcal{F}$ of $\Pi_1$ shares the same line with $\Pi_2$ and $\Pi_3$. Then some points of this line necessarily lie on three elements of $\mathcal{F}$ (false). So $\Pi_{1} \cap
  \Pi_{2} \cap \Pi_{3}$ intersects the solid containing the hyperbolic quadric $\mathcal{L}_{1} \cup \mathcal{L}_{2}$ in a plane. This plane contains already one line $\ell_1$ of this hyperbolic quadric $\mathcal{L}_{1} \cup \mathcal{L}_{2}$, so it contains a second line of $\mathcal{L}_{1} \cup \mathcal{L}_{2}$.
  
   But for each $\pi_{1} \in \Pi_{1}$, we
  find that the line $\pi_{1} \cap \Pi_{2}$ cannot lie in $\Pi_{3}$ since
  otherwise $\pi_{1} \cap \Pi_{2}$ would meet $q$ elements of $\mathcal{F}$
  in $\Pi_{2}$ and $q$ elements of $\mathcal{F}$ in $\Pi_{3}$, a
  contradiction. Thus $\ell_{2} = \bar{\pi}_{2} \cap \Pi_{1}$ must be the second
  line of the hyperbolic quadric in $\Pi_{3}$. 

  By symmetry, we also find that $\bar{\pi}_{3}$ intersects $\Pi_{1}$ and
  $\Pi_{2}$ in the same line.

  Applying this argument for all the
  $\frac{1}{2}[\frac{q^{n}-1}{q-1}-\delta]+1$ different parallel spaces
  found in the first step, we obtain a space
  $\Omega$ in the common intersection which contains
  $\frac{1}{2}[\frac{q^{n}-1}{q-1}-\delta]+1$ different lines of contact
  points. This proves that $\Omega$ must have dimension $n$ and we can copy
  the final steps of the case $q$ odd to prove that $\Omega$ contains only
  contact points. \qed
\end{Proof}

\noindent {\bf Concluding arguments.}
Applying Lemma~\ref{l:fine} precisely $\delta$ times, we find that $\mathcal{F}$ can
be extended to a dual arc $\mathcal{F}'$ of size
$\frac{q^{n+1}-1}{q-1}$. Even in the case $q$ even, no $2n$-dimensional space
contains exactly $2$ elements of $\mathcal{F}'$. By
Result~\ref{JefHendrik}, this implies 
that $\mathcal{F}'$ is the dual arc given by Construction
\ref{cons:1}. As we know from Theorem~\ref{T:d1-ext}, in the case $q$
even this dual arc 
can be extended by one extra element.

\medskip

\noindent
{\bf Acknowledgements} 
The authors profoundly like to thank W.~M.~Kantor and E.~E.~Shult for their helpful comments on an earlier version of this paper. The research of the first and the second author took place within the project "Linear codes and cryptography" of the Fund for Scientific Research Flanders (FWO-Vlaanderen) (Project nr. G.0317.06), and the research of the three authors is supported by the Interuniversitary Attraction Poles Programme-Belgian State-Belgian Science Policy: project P6/26-Bcrypt.\\ 

Address of the authors:\\

\noindent A. Klein and L. Storme: Ghent University, Dept. of  Mathematics, Krijgslaan 281-S22, 9000 Ghent, Belgium\\

\noindent J. Schillewaert: 
D\'epartement de Math\'ematique,
Universit\'e Libre de Bruxelles, 
U.L.B., CP 216,
Bd. du Triomphe,
B-1050 Bruxelles,
Belgique\\

\noindent A. Klein: klein@cage.ugent.be, \url{http://cage.ugent.be/~klein}

\noindent J. Schillewaert: jschille@ulb.ac.be 

\noindent L. Storme: ls@cage.ugent.be, \url{http://cage.ugent.be/~ls}

\end{document}